\documentstyle[amssymb,12pt]{article}

\input amssym.def
\input amssym

\title {\Huge{An introduction to o-minimal structures}
\thanks{To appear in the proceedings of the TMR 
Junior and Summer School in Complex Dynamics, 2-11 September 1999, 
CMAF Lisboa, Portugal.}}

\author {M\'{a}rio J. Edmundo \thanks{Partially supported by JNICT grant 
PRAXIS XXI/BD/5915/95 and a TMR grant for the meeting.}
\\The Mathematical Institute, 
\\24-29 St Giles
\\OX1 3LB
\\Oxford, U.K
\\edmundo@maths.ox.ac.uk\\}
\date{September 3, 1999}

\newcommand{\into}{\longrightarrow}

\renewcommand{\tilde}{\widetilde}
\renewcommand{\bar}{\overline}

\newcommand{\QQ}{\mathbb{Q}}
\newcommand{\RR}{\mathbb{R}}
\newcommand{\CC}{\mathbb{C}}

\newcommand{\Ll}{\mbox{$\mathbf L$}}
\newcommand{\M}{\mbox{$\mathbf M$}}
\newcommand{\N}{\mbox{$\mathbf N$}}
\newcommand{\R}{\mbox{$\mathbf R$}}

\newcommand{\Pp}{\mbox{$\mathbf P$}}

\begin{document}

\maketitle
\begin{abstract}
The first papers on o-minimal structures appeared in the mid
1980s, since then the subject has grown into a wide ranging
generalisation of semialgebraic, subanalytic and subpfaffian
geometry. In these notes we will try to show that this is in fact the
case by presenting several examples of o-minimal structures and by
listing some geometric properties of sets and maps definable in
o-minimal structures. We omit here any reference to the pure model
theory of o-minimal structures and to the theory of groups and rings
definable in o-minimal structures.
\end{abstract}

\newpage
{\small
\begin{section}{Basic model theory}\label{section basic model theory}

We start by recalling  some basic notions of model theory (structures, 
expansions and reducts, definable sets, and elementary extensions, etc., - the
reader already familiarised with these notions can skip this section) and
we end the section by illustrating with examples these model theoretic 
notions.  
\[ \]
\textbf{Basic model theory.}
A {\it structure} $\N$ consists of: $(1)$ a non empty set $N$; $(2)$ a
set of constants $(c_k^N)_{k\in K}$, where $c_k^N\in N$; $(3)$ a family
of maps $(f_j^N)_{j\in J}$, where $f_j^N$ is an $n_j$-ary map,
$f_j^N:N^{n_j}\into N$ and  $(4)$ a family $(R_i^N)_{i\in I}$ of 
relations that is, for each $i$, $R_i^N$ is a subset of
$N^{n_i}$ for some $n_i\geq 1$. We often use the following notation
$\N$$=(N, (c_k^N)_{k\in K},(f_j^N)_{j\in J}, (R_i^N)_{i\in I})$, and 
sometimes we omit the superscripts. The {\it language} $\Ll$
associated to a structure $\N$ consists of:  $(1)$ For each constant
$c_k^N$, a constant symbol $c_k$; $(2)$ For each map $f_j^N$ a
function symbol, $f_j$ of arity $n_j$ and  $(3)$ For each relation
$R_i^N$, a relation symbol, $R_i$ of arity $n_i$. We also
include in $\Ll$ a countable set of variables $(x_q)_{q\in Q}$
and use the notation
$\Ll$$=\{ (c_k)_{k\in K}, (f_j)_{j\in J}, (R_i)_{i\in I}, (x_q)_{q\in Q}\}$. If
$\Ll$ is the language associated with the structure $\N$ we say that $\N$ is an
{\it $\Ll$-structure}. If $\Ll$$'\subseteq \Ll$ are two languages and
$\N '$ and $\N$ are respectively an $\Ll$$'$-structure and an
$\Ll$-structure, such that $N'=N$ then we say that $\N$ is an 
{\it expansion} of $\N '$ or that $\N '$ is a {\it reduct} of $\N$.

Let $\Ll$ be a language and $\N$ an $\Ll$-structure. We are going to
define inductively the set of {\it $\Ll$-formulas} and
{\it satisfaction} of an $\Ll$-formula $\phi $ in the $\Ll$-structure $\N$,
in order to define the {\it $\N$-definable sets}. The set of
{\it $\Ll$-terms} is generated inductively by the following rules: $(i)$
every variable is an $\Ll$-term, $(ii)$ every constant of $\Ll$ is an
$\Ll$-term and $(iii)$ if $f$ is in $\Ll$ is an $n$-ary function, and 
$t_1,\dots ,t_n$ are $\Ll$-terms, then $f(t_1,\dots ,t_n)$ is an
$\Ll$-term. An {\it atomic} $\Ll$-formula is an
expression of the form: $t_1=t_2$ or $R(t_1,\dots ,t_n)$ where $R$ is
an $n$-ary relation in $\Ll$ and $t_1,\dots ,t_n$ are $\Ll$-terms. We 
sometimes write $R(t_1(x_1,\dots ,x_k),\dots ,t_n(x_1,\dots ,x_k))$ if
we want to explicitly show the variables occurring in the atomic
$\Ll$-formula. Given a tuple $a\in N^k$, we say that $a$ {\it
satisfies} $R(t_1(x),\dots ,t_n(x))$ in $\N$ where $x=(x_1,\dots
,x_k)$ if $R^N(t_1(a),\dots ,t_n(a))$ holds. We denote this by 
$\N$$\models R(t_1(a),\dots ,t_n(a))$.

We say that $S\subseteq N^k$ is an {\it atomic $\N$-definable subset 
(defined over $A\subseteq N$)} if there is $b\in A^m$ such that 
$S=\{a\in N^k:\N$$\models \phi (a,b)\}$ for some atomic $\Ll$-formula 
$\phi (x,y)$ with $x=(x_1,\dots ,x_k)$ and $y=(x_{k+1},\dots ,x_{k+m})$ .
We now generate the $\Ll$-formulas (resp., the $\N$-definable sets)
from the atomic $\Ll$-formulas (resp., atomic $\N$-definable sets)
using the following operators: $\wedge $ (and) which corresponds to
intersection, $\vee $ (or) which corresponds to union, $\neg $ (not)
which corresponds to complementation, $\exists $ (there exists)
corresponding to projection and $\forall $ (for all) which corresponds
to inverse image under projections. The construction is as follows:
$(i)$ all atomic $\Ll$-formulas are $\Ll$-formulas; $(ii)$ if $\phi
_1(x)$ and $\phi _2(x)$ are $\Ll$-formulas, then $(\phi _1\wedge \phi
_2)(x)$ and $(\phi _1\vee \phi _2)(x)$  are $\Ll$-formulas, and for
$a\in N^k$, $\N$$ \models (\phi _1\wedge \phi _2)(a)$ iff $\N$$\models
\phi _1(a)$ and $\N$$\models \phi _2(a)$, we also have the obvious
clause for $\vee $; $(iii)$ if $\phi (x)$ is an $\Ll$-formula, $\neg
\phi (x)$ is an $\Ll$-formula, and for $a\in N^k$, $\N$$ \models  \neg
\phi (a)$ iff $\phi (a)$ does not hold in $\N$ (this is denote by $\N$
$\nvDash \phi (a)$); $(iv)$ if $\phi (x,x_{k+1})$ is an $\Ll$-formula,
then $\exists x_{k+1}\phi (x,x_{k+1})$ is an $\Ll$-formula, and for  
$a\in N^k$, $\N$$\models \exists x_{k+1} \phi (a,x_k)$ iff there
exists $b\in N$ such that $\N$$\models \phi (a,b)$; $(v)$ the obvious 
clauses  for $\forall $. The $\Ll$-formulas constructed using only
$(i)$, $(ii)$ and $(iii)$ are called {\it quantifier-free
$\Ll$-formulas}. Two $\Ll$-formulas $\phi _1(x)$ and $\phi _2(x)$
are {\it equivalent} in $\N$ if for all $a\in N^k$ $\N$$ \models \phi _1(a)$ 
iff $\N$$\models \phi _2(a)$. We say that $\N$
has {\it quantifier elimination} if every $\Ll$-formula is equivalent
in $\N$ to a quantifier-free $\Ll$-formula; we say that $\N$ is {\it
model complete} if every $\Ll$-formula is equivalent
in $\N$ to an {\it existential} $\Ll$-formula i.e., an $\Ll$-formula
of the form $\exists y\phi (x,y)$.  
 
A subset $D\subseteq N^k$ is an {\it $\N$-definable subset (defined
over $A\subseteq N$)} if there is an $\Ll$-formula $\phi (x,y)$ with 
$x=(x_1,\dots ,x_k)$ and $y=(x_{k+1},\dots ,x_{k+m})$ and some $b\in
A^m$ such that $D=\{a\in N^k:\N$$\models  \phi (a,b) \}$. The
``$\N$-constructible sets'' are those $\N$-definable sets determined
by some quantifier-free $\Ll$-formula, they are finite boolean
combination of atomic $\N$-definable sets. If $\N$ has quantifier
elimination then every $\N$-definable set is an $\N$-constructible set
and if $\N$ is model complete then every $\N$-definable set is a
projection of an $\N$-constructible set.
If $A\subseteq N^k$ and $B\subseteq N^m$ are $\N$-definable sets (over 
$A\subseteq N$), a function $f:A\into B$ is $\N$-definable (over $A$)
if its graph is an $\N$-definable set (over $A$). More generally, a 
structure $\M$$=(M, (c_k^M)_{k\in K}, (f_j^M)_{j\in J}, (R_i^M)_{i\in
I})$ is $\N$-definable (over $A$) if: $(i)$ $M\subseteq N^l$ is $\N$-
definable (over $A$); $(ii)$ for each $k\in K$ there is  a point
$m_k\in M$ corresponding to $c_k^M$; $(iii)$ for each $j\in J$ the
function $f_j^M:M^{n_j}\into M$ is $\N$-definable (over $A$) and
$(iv)$ for each $i\in I$ the relation $R_i^M \subseteq M^{n_i}$ is
$\N$-definable (over $A$). Note that, in this case every
$\M$-definable set is also an $\N$-definable set. 

Given two $\Ll$-structures $\N$ and $\M$, a map $h:N\into M$ (which
determines in the obvious way a map $h:N^k\into M^k$) is a {\it
homomorphism} if : $(i)$  for every constant $c$ in $\Ll$,
$h(c^N)=c^M$; $(ii)$ for every $n$-ary function $f$ in $\Ll$, for
every $a\in N^n$, $h(f^N(a))=f^M(h(a))$ and $(iii)$ for every $n$-ary 
relation $R$ in $\Ll$, for every $a\in N^n$,  if $R^N(a)$ then
$R^M(h(a))$. An injective homomorphism $h:N\into M$ is an {\it
embedding} if for every $n$-ary relation $R$ in $\Ll$, for every $a\in
N^n$, $R^N(a)$ if and only if $R^M(a)$. An {\it isomorphism} is a
bijective embedding. We say that $\N$ is an {\it $\Ll$-substructure}
of $\M$, denoted by $\N$$\subseteq \M$, if $N\subseteq M$ and the
inclusion map is an embedding. Let $\N$$\subseteq \M$. We say that
$\M$ is an {\it elementary extension} of $\N$ (or that $\N$ is an {\it 
elementary substructure} of $\M$), denoted by $\N$$\preceq \M$, if for
every $\Ll$-formula $\phi (x)$, for all $a\in N^k$, we have
$\N$$\models \phi (a)$ iff  $\M$$\models \phi (a)$. This is equivalent 
(by Tarski-Vaught test) to saying that for every non empty
$\M$-definable set $E\subseteq M^l$, defined with parameters from $N$, 
$E(N):=E\cap N^l$ ("the set of $N$-points of $E$") is a non empty
$\N$-definable set. Clearly, if $S\subseteq N^l$ is an $\N$-definable
set defined with parameters from $N$ and $\N$$\preceq \M$, then the
$\Ll$-formula which determines $S$ determines an $\M$-definable set 
$S(M)\subseteq M^l$ ("the $M$-points of $S$"). 
The {\it theory} $Th(\N$$)$ of an $\Ll$-structure $\N$ is the
collection of all $\Ll$-sentences (i.e., $\Ll$-formulas without free
variables) $\sigma $ such that $\N$$ \models  \sigma $. $\N$ is {\it 
elementarily equivalent} to $\M$, denoted $\N$$\equiv \M$ iff
$Th(\N$$)=Th(\M$$)$. Clearly, if $\N$$\preceq \M$ then $\N$$\equiv
\M$. Note also that if $\N$ has quantifier elimination (resp., is
model complete) and $\N$$\equiv \M$ then $\M$ has quantifier
elimination (resp., is model complete).

The following two  facts (the L\"{o}wenheim-Skolem theorems) are
fundamental theorems of basic model theory: $(1)$ Let $\Ll$ be a
language, $\N$ an $\Ll$-structure and $X\subseteq N$. 
Then for every cardinal $\kappa $ such that $|X|+|\Ll$$|\leq \kappa
\leq |N|$, $\N$ has an elementary substructure $\M$ such that
$X\subseteq M$ and $|M|=\kappa $;
$(2)$ Let $\Ll$ be a language, let $\N$ be an infinite
$\Ll$-structure. Then for any cardinal $\kappa >|N|$, $\N$ has an
elementary extension of cardinality $\kappa $. 

A set of $\Ll$-sentences $\Sigma $ is {\it consistent} if there is an 
$\Ll$-structure $\N$ such that
for all $\sigma \in \Sigma $ we have $\N$$\models \sigma $. In this
case we say that $\N$ is a {\it model} of $\Sigma $, denoted
$\N$$\models \Sigma $. An {\it $\Ll$-theory} (resp., a {\it complete
$\Ll$-theory}) is a consistent set of $\Ll$-sentences (resp., a maximal  
consistent set of $\Ll$-sentences). An $\Ll$-theory $T$ is {\it
axiomatizable} if there is a set of $\Ll$-sentences $\Sigma $ (the set
of axioms) such that for every $\Ll$-structure $\N$, $\N$$\models T$
iff $\N$$\models \Sigma $.
The {\it Compactness theorem} says that: if $\Sigma $ is a set of 
$\Ll$-sentences then $\Sigma $ is consistent iff every finite subset
of $\Sigma $ is consistent. Moreover, if $\N$ is an $L$-structure and
$F$ is a family of $\N$-definable subsets of $N^k$ with finite
intersection property in $N^k$ then there is an elementary
extension $\M$ of $\N$ such that $F$ has non empty intersection in
$M^k$.
\[ \]
\textbf{Examples.}
Let $\Ll$$_{rings}:=\{ 0,1, +, -,\cdot \}$ be the language of
rings. Then $\CC$$:=(\CC$$, 0,1,+$$, -,\cdot )$ is an
$\Ll$$_{rings}$-structure. The atomic sets in $\CC$$^l$ are exactly the 
Zariski closed sets, and by Chevalley's theorem, the $\CC$-definable sets are
the constructible sets (i.e., boolean combinations of Zariski closed
sets) which means that $\CC$ has quantifier elimination. The theory
$Th(\CC$$)$ of
$\CC$ is the theory $ACF_0$ of {\it algebraically closed fields}
of characteristic zero
which is axiomatised by the usual axioms for fields of characteristic
zero, together with 
$\forall x_1\dots \forall x_l\exists y (y^l+x_1y^{l-1}+\cdots +
x_{l-1}y+x_l=0)$, for each positive integer $l$. Another model of
$ACF_0$ is the algebraic closure of $\QQ$ in $\CC$. The compactness 
theorem shows that $ACF_0$ has model of any transcendence degree and
the L\"{o}weinheim-Skolem theorem shows that $ACF_0$ has models of
any infinite cardinality. Two model of $ACF_0$ are isomorphic iff
their transcendence base over $\QQ$ has the same cardinality, in
particular for an uncountable cardinal $\kappa $, up to isomorphism
there is only one model of $ACF_0$ of cardinality $\kappa $. (And
there are $2^{\aleph _0}$ countable models).
The models of $ACF_0$ are called {\it algebraically closed fields} 
(of characteristic zero) and are examples of  {\it strongly
minimal structures} i.e., structures $\N$ such that any $\N$-definable
subset of $N$ is either finite or co-finite. 
Other examples of such structures are a  nonempty set and a vector
space over a division ring.

Let $\Ll$$_{ord}:=\{ 0,1,+, -,\cdot ,< \}$ be the language of ordered
rings. Then $\RR$$:=(\RR$$, 0,1,+$$, -,\cdot ,< )$ is an
$\Ll$$_{ord}$-structure. Tarski \cite{T} showed that $\RR$ has
quantifier elimination and so the $\RR$-definable subsets of $\RR$$^l$
are boolean combinations of sets of the form 
$\{\bar{a}\in \RR$$^l:f(\bar{a})=0\}$ and 
$\{\bar{a}\in \RR$$^l:g(\bar{a})>0\}$ where 
$f,g\in \RR$$[x_1,\dots ,x_l]$. These sets are called semi-algebraic. 
The theory $Th(\RR$$)$ of $\RR$ is the theory $RCF$ of {\it real closed
fields} which is axiomatised by the usual axioms for ordered 
fields together with (i) $\forall x_1\dots \forall x_l (x_1^2+\cdots +
x_l^2\neq 1)$ (for each positive integer $l$), (ii) $\forall x\exists
y (x=y^2 \vee -x=y^2)$ and (iii) $\forall x_1\dots \forall x_l\exists
y (y^l+x_1y^{l-1}+\cdots +x_{l-1}y+x_l=0)$ (for all odd $l$). 
Another model of $RCF$ is the algebraic closure of $\QQ$ in $\RR$.
Again, the L\"{o}weinheim-Skolem theorem shows that $RCF$ has models of
any infinite cardinality, but unlike $ACF$, for any infinite
cardinal $\kappa $, up to isomorphism there are $2^{\kappa }$ model of 
$RCF$ of cardinality $\kappa $.
The models of $RCF$ are called {\it real closed fields} 
and they are examples of o-minimal structures (see the definition below). 
\end{section}

\begin{section}{O-minimal structures}\label{section o-minimal structures}

\textbf{O-minimal structures.}
An {\it o-minimal structure} is an expansion $\N$$=(N,<,\dots )$ of a 
linearly ordered nonempty set $(N,<)$, such that every $\N$-definable
subset of $N$ is a finite union of points and intervals with endpoints in 
$N\cup \{-\infty, +\infty \}$. 

Note the following important results:
let $\N$ be an o-minimal structures then: $(0)$ every $\N$-definable
structure $\M$ which is an expansion of a linearly ordered nonempty set
$(M,<_M)$ is also o-minimal; $(1)$ \cite{KPS} if $\M$ is a structure (in the
language of $\N$) such that $\N$$\equiv \M$ then $\M$ is also o-minimal; $(2)$
\cite{PiS1} for every $A\subseteq N$ there is a {\it prime model} of $Th(\N$$)$
over $A$ (or simply prime model if $A$ is empty) i.e., there is an o-minimal
structure $\Pp$ such that $A\subseteq P$, $\Pp$$\equiv \N$, $\Pp$ is unique up
isomorphism over $A$ and for all $\M$$\equiv \N$ with $A\subseteq M$,
there is an elementary embedding $P\into M$ which is the identity over $A$. In
particular by $(2)$, if $S$ is an $\N$-definable set over $P$ and $\M$
$\equiv \N$ with $A\subseteq M$, then $S$ determines an $\M$-definable
set $S(M)$. 
Let $\N$ be an o-minimal structure in the language $\Ll$. For every 
$\kappa >\max\{\aleph _0, |\Ll$$|\}$ there up to isomorphism
$2^{\kappa }$ o-minimal structures $\M$ such that $|M|=\kappa $ and 
$\M$$\equiv \N$ (see \cite{Sh}), and
if $\Ll$ is countable then up to isomorphism there are either
$2^{\aleph _0}$ or
$6^n3^m$ countable o-minimal structures $\M$ such that $\M$$\equiv \N$
(see \cite{M}).

Let $\N$ be an o-minimal structure. We now list some geometric
properties of $\N$-definable sets and $\N$-definable maps, most of
these can be found in book \cite{D2} but the proof appear elsewhere. 
Note that because of the presence of an ordering in $N$, each
$N^m$ has a natural topology and if $\N$ is an expansion of ordered
ring it makes sense to talk about differentiability (but in general,
it does not make sense to talk about integrability).  Two of the most 
powerful results are the $C^p$-cell decomposition theorem for
$\N$-definable sets and $\N$-definable maps (where $p=0$ if $\N$ is
not an expansion of an ordered ring) and the monotonicity theorem for 
$\N$-definable one variable functions. The cell decomposition theorem
has several consequences: (1) it is  used to define the notion of {\it 
o-minimal dimension} and {\it o-minimal Euler characteristic} for
$\N$-definable sets, these notions are well behaved under the usual
set theoretic operations on $\N$-definable sets, are invariant under
$\N$-definable bijections and given an $\N$-definable family of 
$\N$-definable sets, the set of parameters whose fibre in the family
has a fixed dimension (resp., Euler characteristic) is also an
$\N$-definable set; (2) it shows that every $\N$-definable set has
only finitely many $\N$-definably connected components, and given an 
$\N$-definable family of $\N$-definable sets there is a uniform bound
on the number of $\N$-definably connected components of the fibres in
the family. 

A local version of the o-minimal analog of Zilber's conjecture holds 
\cite{PeS}: if $a\in N$ then the structure induced by $\N$ on an open 
interval containing $a$ is either trivial, or the structure of an open 
interval in an ordered vector space over some ordered division ring or
an o-minimal expansion of an ordered ring. If $\N$ expands an ordered
group, then $\N$ is either {\it linear}, {\it eventually linear}, or
{\it linearly bounded} (several useful characterisations of these
three cases are given in \cite{LP}, \cite{E1} and \cite{MS}).  When
$\N$ expands an ordered ring then $\N$ is either {\it power bounded}
({\it polynomially bounded} if $N$ is Archemedian) or is {\it
exponential} (see\cite{Mi1}) moreover, several geometric 
properties from semialgebraic and subanalytic geometry also
hold for $\N$-definable sets and $\N$-definable maps: we have (1)  an 
$\N$-definable
curve selection theorem (this in fact holds in the more general case
where $\N$ has $\N$-definable Skolem functions e.g., $\N$ expands an
ordered group); (2) an $\N$-definable triangulation theorem; (3) an 
$\N$-definable trivialization theorem; (4) and finally (see
\cite{DM2}) an $\N$-definable analog of the uniform bounds on growths 
theorem, the $C^p$-multiplier theorem, the generalised Lojasiewicz 
inequality, the $C^p$ zero set theorem, the $C^p$ Whitney
stratification theorem, etc. 
\[ \]
\textbf{Examples.}
As examples of o-minimal structures apart from the trivial ones such
as: (1) dense linearly ordered nonempty sets without endpoints
-conversely by \cite{PiS1}, if $\N$$=(N,<,\dots )$ is an o-minimal
structure, then $(N,<)$ is elementarily equivalent to an ordered set
of the form $C_1+\cdots +C_m$ where: (i) $C_i$ is elementarily
equivalent to one of the following ordered sets: a finite ordered set, 
$\omega $, $\omega ^{*}$, $\omega +\omega ^{*}$, $\omega ^{*}+\omega
$, $\QQ$ and (ii) if $C_i$ does not have a last element, then
$C_{i+1}$ has a first element, also by results from \cite{PiS2} and 
\cite{PiS3} one usually assumes, without loss of generality, that
$(N,<)$ is a dense linearly ordered set without end points; (2) ordered
divisible abelian groups, in fact also ordered vector spaces over a division
ring (semilinear geometry) - and conversely by \cite{PiS1} if $\N$ is an
o-minimal expansion of an ordered group $(N,0,+,<)$ then $(N,0,+,<)$ is an
ordered divisible abelian group, in fact it is also an ordered vector spaces
over the ordered division ring $\Lambda (\N$$)$ of all $\N$-definable 
endomorphisms of 
$(N,0,+,<)$ - we also have the following examples which are of special 
interest to geometers: 

(3) $\RR$$:=(\RR$$,0,1,+,\cdot ,<)$ (semialgebraic 
geometry, by Tarski-Seidenberg theorem \cite{T} this structure has quantifier 
elimination and therefore every $\RR$-definable set is a semialgebraic
set and o-minimality follows from this, in fact any real closed field is 
o-minimal i.e., any ordered ring $\R$$:=(R,0,1,+,\cdot ,<)$ such that 
$\R$$\equiv \RR$, for example the algebraic closure of $\QQ$ in $\RR$
is a real closed field - and conversely by \cite{PiS1} if $\N$ is an
o-minimal expansion of an ordered ring $(N,0,1,+,\cdot ,<)$ then 
$(N,0,1,+,\cdot ,<)$ is a real closed field; 

(4) $\RR$$_{an}:=(\RR$
$,0,1,+,\cdot ,<,(f)_{f\in an})$ where $an$ is the collection of all 
functions which are the restriction to $[-1,1]^n$ (for some $n$) of analytic 
functions on some open neighbourhood of $[-1,1]^n$ (by results of Gabrielov 
\cite{G}, Lojasiewicz \cite{L} and Bierstone and Milman \cite{BM} the 
$\RR$$_{an}$-definable sets are exactly the subanalytic sets in the projective 
spaces we therefore get (global) subanalytic geometry, this structure is 
model complete and o-minimal as remarked by van den Dries \cite{D1};
$\RR$$_{an}$ has quantifier elimination after adding the function
$1/x$ (with $1/0=0)$ to its language, an axiomatization of its theory
is given in \cite{DMM1}. 
Some model complete reducts of $\RR$$_{an}$  expanding $\RR$ were 
constructed: the expansions of $\RR$ by restricted elementary
functions \cite{D3} and expansions of $\RR$ by abelian and elliptic 
functions. By \cite{Mi2}, the expansion of $\RR$$_{an}$ by power
functions has quantifier elimination, is o-minimal and has an explicit
axiomatization;   

(5) Wilkie \cite{W1} uses model theory, valuation theory 
and results by Khovanskii \cite{K} to show that $\RR$$_{exp}:=(\RR$
$,0,1,+,\cdot ,<,exp )$  and the expansion of $\RR$ by restricted Pfaffian 
functions are model complete and o-minimal;  By \cite{Mi2}, the
expansion of $\RR$ by restricted Pfaffian functions, power functions
and a constant symbol for each exponent is model complete and o-minimal.
A quantifier elimination result and a (non trivial) axiomatisation of
the expansion of $\RR$ by the restricted Pfaffian functions (and also
by power functions and a constant symbol for each exponent) is not
known. Ressayre gives an axiomatisation of the theory of
$\RR$$_{exp}$, as for quantifier elimination, van den Dries \cite{D2}
adapts an old result of Osgood to show that an expansion of $\RR$ by a 
family of total real analytic functions admits elimination of
quantifiers iff each such function is semialgebraic.  Macintyre and
Wilkie \cite{MW} show that if the Schanuel conjecture holds then 
$Th(\RR$$_{exp})$ is decidable;

(6) Wilkie's method and 
Khovanskii result are refined in \cite{DM1} to show that the structure
$\RR$$_{an,exp}
:= (\RR$$,0,1,+,\cdot ,<, (f)_{f\in an}, exp)$ 
is model complete and o-minimal, van den Dries, Macintyre and Marker 
\cite{DMM1}, inspired by work of 
Ressayre (a preliminary version of \cite{Re}) give a different proof of this 
fact and in \cite{DMM2} explicit nonstandard models of $\RR$$_{an,exp}$ are 
constructed leading to a solution of a conjecture posed by Hardy and to some 
non definability results: the following functions are not $\RR$$_{an,exp}$-
definable: $\Gamma _{|(0,+\infty)}$ (the gamma function), the error function, 
the logarithmic integral and $\zeta _{|(1,+\infty )}$ (the Riemann zeta 
function). By \cite{DMM1}, the expansion $\RR$$_{an,exp,log}$ of 
$\RR$$_{an,exp}$ has quantifier elimination and both have an explicit 
axiomatisation.
A geometric proof of o-minimality and model completeness of 
$\RR$$_{an,exp}$ has been given recently by Lion and Rolin \cite{LR1};

(7) Denef and van den Dries \cite{DD} give a proof of model completeness 
and o-minimality of $\RR$$_{an}$ using more explicitly the Weirstrass 
preparation theorem, this is then generalised to establish the model 
completeness and o-minimality of $\RR$$_{an^*}$ the expansion of $\RR$ 
by (restricted convergent) generalized power series (\cite{DS1}) and
of $\RR$$_{\mathcal{G}}$ the 
expansion of $\RR$ by (a variant of) Tougeron's class of Gevrey functions 
(\cite{DS2}). In \cite{DS2} the expansions $\RR$$_{an^*,exp}$ and 
$\RR$$_{\mathcal{G}}$$_{,exp}$ of $\RR$$_{an^*}$ and $\RR$$_{\mathcal{G}}$ by 
$exp$ are shown to be o-minimal and model complete, in particular 
$\zeta _{|(1,+\infty )}$ is definable in $\RR$$_{an^*,exp}$ since 
$\zeta (-log(x))=\sum _{n=1}x^{logn}$ is definable in $\RR$$_{an^*}$ and 
$\Gamma _{|(0,+\infty )}$ is definable in $\RR$$_{\mathcal{G}}$$_{,exp}$
since $log\Gamma (x)=(x-\frac{1}{2})logx-x+\frac{1}{2}log(2\pi )+\phi (x)$
where $\phi $ is definable in $\RR$$_{\mathcal{G}}$; Quantifier
elimination results and (non trivial) axiomatisations for
$\RR$$_{an^*}$ and $\RR$$_{\mathcal{G}}$ are not known. But by
\cite{DS2}, if $\tilde{\RR}$ is a polynomially bounded o-minimal
expansion of $\RR$ such that $exp_{|[0,1]}$ is
$\tilde{\RR}$-definable, then the expansion of $\tilde{\RR}$ by $exp$
(resp., $exp$ and $log$) is model complete (resp., has quantifier
elimination) and is o-minimal, moreover they have explicit
axiomatisations - the axiomatisation of $Th(\tilde{\RR}$$)$ plus
Ressayre axioms for $exp$ (resp., $exp$ and $log$)-
and they are exponentially bounded. 

(8) Finally, building
on work of Charbonnel, Wilkie \cite{W2} gives necessary and 
sufficient conditions for an expansion of $\RR$ by total $C^{\infty }$ 
functions to be o-minimal, in particular
o-minimality of the expansion of $\RR$ by total 
$C^{\infty }$ Pfaffian functions is established, a geometric treatment
of Wilkie's result in the subanalytic context is given by Lion and
Rolin \cite{LR2} using Moussu and Roche's \cite{MR} notion of Rolle
leafs and the Khovanskii-Rolle theorem, 
this is later generalised \cite{Sp} to show that any o-minimal expansion 
$\tilde{\RR}$ of $\RR$ has an o-minimal Pfaffian closure $\mathcal{P}$
$(\tilde{\RR}$$)$ i.e., the Rolle leafs of $1$-forms with $C^1$ 
coefficients definable in $\mathcal{P}$$(\tilde{\RR}$$)$ are already 
definable in $\mathcal{P}$$(\tilde{\RR}$$)$. In \cite{LS} the prove of
the existence of the o-minimal Pfaffian closure
$\mathcal{P}$$(\tilde{\RR})$ of $\tilde{\RR}$ is refined to show that
if $\tilde{\RR}$ has analytic cell decomposition ( resp., is
exponentially bounded) then $\mathcal{P}$$(\tilde{\RR})$
has analytic cell decomposition (resp., is exponentially bounded). Its
not known if the model completeness of $\tilde{\RR}$ implies that of 
$\mathcal{P}$$(\tilde{\RR})$, but a relative model completeness result
is proved in \cite{LS}.

\medskip
All the examples of o-minimal expansions of $\RR$ mentioned above have 
analytic cell decomposition and those which do not have $exp $ on
their language are polynomially bounded (the other ones are
exponentially bounded). There is no known example of an exponential o-minimal
expansion of $\RR$ which is not exponentially bounded.

\end{section}

}

\begin{thebibliography}{XXXXXXX}

\bibitem[BM]{BM}E.Bierstone and P.Milman {\it Semi-analytic and
subanalytic sets} IHES Publ. Math. 67 (1988) 5-42.

\bibitem[D1]{D1}L. van den Dries {\it A generalisation of the 
Tarski-Seidenberg theorem} Bulletin of the AMS 15 (1986) 189-193.

\bibitem[D2]{D2}L. van den Dries {\it Tame Topology and o-minimal 
structures} Cambridge University Press 1998.

\bibitem[D3]{D3}L. van den Dries {\it On the elementary theory of
restricted elementary functions} JSL 53 (1988) 796-808.

\bibitem[DD]{DD}J.Denef and L.van den Dries {\it $p$-adic and real
subanalytic sets} Annals of Math. 128 (1988) 79-138.

\bibitem[DM1]{DM1}L. van den Dries and C.Miller {\it On the real
exponential field with restricted analytic functions} Israel
J. Math. 85 (1994) 19-56.

\bibitem[DM2]{DM2}L. van den Dries and C.Miller {\it Geometric
categories and o-minimal structures} Duke Math. Journal (2) 82 (1996) 497-540.

\bibitem[DMM1]{DMM1}L. van den Dries, A. Macintyre and D. Marker
{\it The elementary theory of restricted analytic fields with
exponentiation} Annals of Math. 140 (1994) 183-205.

\bibitem[DMM2]{DMM2}L. van den Dries, A. Macintyre and D. Marker
{\it Logarithmic-exponential power series} Journal of the LMS (2) 56
(1997) 417-434.

\bibitem[DS1]{DS1}L. van den Dries and P. Speissegger {\it The
real field with convergent generalised power series} 
Trans. of the AMS (11) 350 (1998) 4377-4421. 

\bibitem[DS2]{DS2}L. van den Dries and P. Speissegger {\it The
field of reals with multisummable series and the exponential function} 
to appear in the Proc. of the LMS (2000+).

\bibitem[E1]{E1}M.Edmundo {\it Structure theorems for o-minimal
expansions of groups} Annals of Pure and Applied Logic 102 (2000) 159-181.

\bibitem[G]{G}A.Gabrielov {\it Projection of semi-analytic sets}
Functional Analysis and its Applications 2 (1968) 282-291.

\bibitem[K]{K}A. Khovanskii {\it On a class of systems of
transcendental equations} Soviet Math. Doklady 22 (1980) 762-765.

\bibitem[KPS]{KPS}J.Knight,A.Pillay and C.Steinhorn {\it Definable
sets in ordered structures.II}
Trans. of the AMS 295 (1986) 593-605.

\bibitem[LR1]{LR1}J.-M.Lion and J.-P.Rolin {\it Th\'{e}or\`{e}me de
pr\'{e}paration pour les fonctions logarithmico-exponentielles} Ann.
Inst. Fourier, 47 (1997) no 3, 859-884.

\bibitem[LR2]{LR2}J.-M.Lion and J.-P.Rolin {\it Volumes, feuilles de
Rolle et feuilletages analytiques reeles et theoreme de Wilkie}
Ann. Toulouse 7 (1998) 93-112.

\bibitem[LS]{LS}J.-M.Lion and P.Speissegger {\it Analytic
stratification in the Pfaffian closure of an o-minimal structure} to
in the Duke Math. Journal (2000+).

\bibitem[L]{L}S.Lojasiewicz {\it Ensembles semi-analytique} Lecture
Notes Ecole Polytechnique, Paris 1965.

\bibitem[LP]{LP}J.Loveys and Y.Peterzil {\it Linear o-minimal
structures} Israel Journal of Math. 81 (1993) 1-30.

\bibitem[MW]{MW}A.Macintyre and A.Wilkie {\it Schanuel's conjecture
implies the decidability of real exponentiation} Preprint 1994.

\bibitem[M]{M}L. Mayer {\it Vaughts conjecture for o-minimal theories} 
JSL 53 (1988) 146-159.

\bibitem[Mi1]{Mi1}C. Miller {\it A growth dichotomy for o-minimal
expansions of ordered fields} Preprint 1996.

\bibitem[Mi2]{Mi2}C. Miller {\it Expansions of the real field with
power functions} Annals of Pure and Applied Logic 68 (1994) 79-94.

\bibitem[MS]{MS}C.Miller and S.Starchenko
{\it A growth dichotomy for o-minimal expansions of ordered groups}
in Logic Colloquium 93, W. Hodges et al., (eds.), Oxford Univ. Press 1994.

\bibitem[MR]{MR} R.Moussu and C.Roche {\it Theorie the Khovanskii et
probleme de Dulac} Invent. Math. 105 (1991) 431-441.

\bibitem[PeS]{PeS}Y.Peterzil and S.Starchenko
{\it A trichotomy theorem for o-minimal structures}
Vanderbilt University Preprint Series No.95-018.

\bibitem[PiS1]{PiS1}A.Pillay and C.Steinhorn {\it Definable sets in 
ordered structures.I} Trans. of the AMS 295 (1986) 565-592.

\bibitem[PiS2]{PiS2}A.Pillay and C.Steinhorn {\it Discrete 
o-minimal structures} Annals of Pure and Applied Logic 34 (1987) 275-290.

\bibitem[PiS3]{PiS3}A.Pillay and C.Steinhorn {\it Definable sets in 
ordered structures.III} Trans. of the AMS 309 (1988) 565-592.

\bibitem[Re]{Re}J.-P. Ressayre {\it Integral part of real closed
exponential fields} in ``Arithmetic, proof theory and computational 
complexity'', Oxford Logic Guides 23, Oxford Univ. Press 1993.

\bibitem[Sh]{Sh}S.Shelah {\it Classification theory and the number of
non isomorphic models} 2nd ed. revised, North-Holland 1990.

\bibitem[Sp]{Sp}P. Speissegger {\it The Pfaffian closure of an
o-minimal structure} J. Reine Angew. Math. 508 (1998) 198-211.

\bibitem[T]{T}A. Tarski {\it A decision method for elementary algebra
and geometry} 2nd ed. revised, Berkeley and Los Angeles 1951.

\bibitem[W1]{W1}A.Wilkie {\it Model completeness results for expansions 
of the ordered field of real numbers by restricted Pfaffian functions 
and the exponential function} Journal of the AMS 9 (1996) 1051-1094.

\bibitem[W2]{W2}A.Wilkie {\it A general theorem of the complement and
some new o-minimal structures} Sel. Math. New ser. 5 (1999) 397-421. 

\end{thebibliography}
\end{document}